\documentclass[a4paper,preprint,10pt]{elsarticle}
\usepackage[english]{babel} 
\usepackage[latin1]{inputenc}
\usepackage{amssymb} 
\usepackage{amsmath}     
\usepackage{pb-diagram}     
\usepackage{epstopdf}                     
\usepackage{graphicx}   
 \usepackage{epsfig}           
\usepackage{color}
\usepackage{fancyhdr}  
 \usepackage{natbib}  
\usepackage{rotating}
\usepackage[T1]{fontenc}        
\usepackage{url}                
\usepackage{mathtools}                 
\usepackage{vmargin}    
\usepackage[title]{appendix}   

\usepackage{lineno}
\usepackage{setspace}
\usepackage{color,soul}
\setulcolor{red} 
\setstcolor{blue}   
\sethlcolor{red}  
\biboptions{square}   
\usepackage[colorlinks=true,linkcolor= blue,urlcolor=blue,citecolor=blue]{hyperref} 
\usepackage[ruled,vlined]{algorithm2e} 
   
\newcommand{\bu}{\bullet}     
\newcommand{\bo}[1]{\mathbf{#1}} 
\def\indic{\hbox{1\kern-.24em\hbox{I}}}      

\newcommand{\x}{x}
\newcommand{\X}{X}

\newcommand{\Z}{Z} 
\newcommand{\R}{\mathbb{R}}

\newcommand{\M}{f}

\newcommand{\T}{T}    
      
\newcommand{\NN}{n}     
  
\newtheorem{prop}{Proposition}{\bf}{\it} 
\newtheorem{defi}{Definition}{\bf}{\it}  

\newtheorem{theorem}{Theorem}{\bf}{\it}     
\newtheorem{lemma}{Lemma}{\bf}{\it}          
\newtheorem{rem}{Remark}{\bf}{\it}  
\newtheorem{corollary}{Corollary}{\bf}{\it} 
\sethlcolor{green}    

\catcode`\"=\active 
\catcode`\"=\active  
\def"{\og\ignorespaces}
\def"{{\fg}} 
 
\makeatletter
\def\ps@pprintTitle{%
  \let\@oddhead\@empty
  \let\@evenhead\@empty
  \def\@oddfoot{\reset@font\hfil\thepage\hfil}
  \let\@evenfoot\@oddfoot
}
\makeatother      
   
\begin{document}

\begin{frontmatter} 
   
  
\title{Gradient and  Hessian of functions with non-independent variables} 
\author[a,b]{Matieyendou Lamboni\footnote{Corresponding author: matieyendou.lamboni[at]gmail.com, matieyendou.lamboni[at]univ-guyane.fr, April 10, 2023}}        
\address[a]{University of Guyane, Department DFRST, 97346 Cayenne, French Guiana, France}
\address[b]{228-UMR Espace-Dev: University of Guyane, University of R\'eunion, IRD, University of Montpellier, France.} 
                                    
\begin{abstract}        
Mathematical models are sometime given as functions of independent input variables and equations or inequations connecting the input variables. A probabilistic characterization of such models results in treating them as functions with non-independent variables. Using the distribution function or copula of such variables that comply with such equations or inequations, we derive two types of partial derivatives of functions with non-independent variables (i.e., actual and dependent derivatives) and argue in favor of the latter. The dependent partial derivatives of functions with non-independent variables rely on the dependent Jacobian matrix of dependent variables, which is also used to define a tensor metric. The differential geometric framework allows for deriving the gradient, Hessian and Taylor-type expansion of functions with non-independent variables.   
\end{abstract}     
                         
\begin{keyword}                         
Copulas \sep Dependent and correlated variables \sep Gradient \sep Hessian \sep Riemannian Tensors
\end{keyword} 
 
               
\end{frontmatter}         
\setpagewiselinenumbers      
  
\section{Introduction} 
We use to work with models defined through functions that include non-independent variables such as correlated input variables. It is also the case of models defined via functions with independent input variables and equations or inequations connecting such inputs; functions subject to some constraints involving input variables and/or the model output. Knowing the key role of and partial derivatives and gradient in mathematical analysis of functions and convergence (\cite{chen15,wang22}), Poincar\'e inequalities and equalities (\cite{bobkov99,fougeres05,bakry08,bobkov09a,bobkov09,roustant17,ren20,lamboni21,lamboni22}), optimization and active subspaces (\cite{russi10,constantine14,zahm20}) and implicit functions (\cite{zhang06,cristea17,li22}), it is relevance to provide i) general-accurate derivatives formulas and gradient of functions which inputs comply with the constraints imposed in particular, and account for the dependence structures among each other in general, ii) the global and probabilistic implicit functions.\\  
                                     
To illustrate one motivation of this work, let us consider the function $\M: \R^3 \to \R$ given by $\M(x, y, z) =2x y^2 z^3$ under the constrained equation $h(x, y, z) =0$ with $h$ any smooth function. Using $\M'_x$ for the formal partial derivative of $\M$ w.r.t. $x$, that is, the partial derivative of $\M$ by considering other inputs as constant or independent, the chain rule gives the following partial derivative:         
\begin{eqnarray}    
\frac{\partial \M}{\partial x}(x, y, z) &=& \M'_x(x, y, z) + \M'_y(x, y, z) \frac{\partial y}{\partial x} 
+ \M'_z(x, y, z) \frac{\partial z}{\partial x} \nonumber \\   
 &=&  2 y^2 z^3 + 4x y z^3 \frac{\partial y}{\partial x}
+ 6x y^2 z^2 \frac{\partial z}{\partial x} \nonumber\, .         
\end{eqnarray}  
 Recall that $(x, y, y)$ are realizations or sample values of a given implicit random  vector $(X, Y, Z)$. So far, the quantities $\frac{\partial y}{\partial x},\, \frac{\partial z}{\partial x}$ and $\frac{\partial \M}{\partial x}$ can be generally hard to determine  without supplementary assumptions. When $(X,\, Z)$ or $(X, \,  Y)$ are independent and using the equation $h(x, y, z) =x +2y -z =0$, we can write 
\begin{equation}            
\displaystyle    
\frac{\partial \M}{\partial x}(x, y, z) =\left\{ \begin{array}{ll}
 2 y^2 z^3 -2 x y z^3  &  \mbox{if $X$ and $Z$ are independent or $z$ is being held fixed}  \\
 2 y^2 z^3 + 6x y^2 z^2 &  \mbox{if $X$ and $Y$ are independent or $y$ is being held fixed}  \\ 
\end{array} \right. \, .        \nonumber     
\end{equation}             
       
The above partial derivatives take account of the relationship between the variables and are also known as actual partial derivatives. It is clear that the actual partial derivative is not unique. While each supplementary assumption can make sense in some cases, it cannot be always guaranteed by the constrained equation $h(x, y, z) =0$ in general. Indeed, when all the initial input variables are dependent or correlated, the above partial derivatives are no longer valid even for a linear function $h$, and it is worth finding the right relationship between the input variables.\\   
        
For non-independent variables having $F$ as the joint cumulative distribution function (CDF), the bi-variate dependency models (\cite{skorohod76}) and  the multivariate dependency models (\cite{lamboni21,lamboni21ar,lamboni22b}) establish formal and analytical relationships among such variables using either CDFs or the corresponding copulas or  new distributions that look like and behave as a copula (\cite{lamboni22mcap}). A dependency function characterizes the probabilistic dependencies among such variables. For a $d$-dimensional random vector of non-independent variables, the dependency models express a subset of $d-1$ variables as a function of independent variables, consisting of the remaining input and $d-1$ new independent variables. Using such dependency models, the contributions of this paper are threefold:  
      
\begin{itemize}           
\item provide a generalization of the actual approach for calculating partial derivatives of function with non-independent variables and give the limits of that approach; 
\item introduce a general-accurate derivative formulas of function with non-independent variables (called dependent partial derivatives);   
\item provide the gradient, Hessian and the Taylor-type expansion of functions with non-independent variables that comply with the dependency structures among the input variables.  
\end{itemize}       
                           
In Section \ref{sec:pcfdv}, firstly, we recall dependency models of dependent input variables, including correlated variables. Secondly, we derive interesting properties of such models regarding the calculus of partial derivatives, and introduce the probabilistic implicit function. By coupling dependency functions with the function of interest, we extend the actual partial derivatives of functions with non-independent variables in Section \ref{sec:fwni}. To avoid the drawbacks of that approach, the dependent partial derivatives of functions with non-independent variables are provided in Section \ref{sec:dgrad}. The gradient and the Hessian matrix of such functions are derived in Section \ref{sec:ugrad} by making use of the differential geometric framework. We  conclude this work in Section~\ref{sec:con}. 


\section*{General notation}                      
For integer $d >0$, let $\bo{\X} :=(X_1, \, \ldots,\, X_d)$ be a random vector of continuous  variables having $F$ as the joint cumulative distribution function (CDF) (i.e., $\bo{\X}  \sim F$). For any $j \in \{1,\ldots, d\}$, we use $F_{x_j}$ or $F_j$ for the marginal CDF of $\X_j$ and  $F_j^{-1}$ for its  inverse. Also, we use $(\sim j) :=(w_1, \ldots, w_{d-1})$ for an arbitrary permutation of $\{1,\ldots, d\}\backslash \{j\}$ and $\bo{\X}_{\sim j} :=(\X_{w_1}, \ldots, X_{w_{d-1}})$. \\      
For a function $\M$ that includes $\bo{\X}$ as inputs, we use $\M'_{x_j}$ for the formal partial derivative of $\M$ w.r.t. $\X_j$, that is, by considering other inputs as constant or independent of $\X_j$, and $\nabla \M := \M'_{\bo{x}} :=\left[\M'_{x_1}, \ldots, \M'_{x_d} \right]^\T$. We use $\frac{\partial \M}{\partial x_j}$ for the partial derivative of $\M$ w.r.t. $\X_j$ that accounts for the dependencies among inputs. We also use $\frac{\partial \M}{\partial \bo{x}} :=\left[\frac{\partial \M}{\partial x_1}, \ldots, \frac{\partial \M}{\partial x_d} \right]^\T \in \R^d$. Of course, we have $\frac{\partial \M}{\partial \bo{x}} = \nabla \M$ for independent inputs. 
                                       
\section{Probabilistic characterization of functions with non-independent variables}\label{sec:pcfdv}   
In probability theory, it is common to treat input variables as random vectors following some CDFs. For instance, for the inputs that take their values within known domains, the Bayesian framework allows for assigning a joint distribution known as \textit{a prior} distribution to such inputs. Without additional information about the inputs, it is common to use  non-informative \textit{a prior} distributions such as uniform distributions or Gaussian distributions with a higher value of the covariance matrix (see e.g., \cite{robert07}). \\ 
     
Functions with non-independent variables include many types of models encountered in practice. An example is the models defined via a given function and equations or inequations connecting its inputs. The resulting inputs that comply with such constraints are often dependent or correlated. In what follows, we are going to use the probability theory for characterizing such non-independent variables. (see Definition \ref{def:fniid}).   
      
\begin{defi} \label{def:fniid}
Consider a function $\M :\R^d \to \R^\NN$ that includes $\bo{\X} \sim F$ as inputs.\\
 
Then, $\M$ is said to be a function with non-independent variables whenever there exists a least a pair $j_1, j_2 \in\{1, \ldots, d\}$ with $j_1 \neq j_2$ such that 
$$
F_{j_1,j_2}(x_{j_1},\, x_{j_2}) \neq F_{j_1}(x_{j_1}) F_{j_2}(x_{j_2}) \, . 
$$   
\end{defi}        
 
We can check that a function that includes $\bo{\X} \sim \mathcal{N}\left(\bo{0}, \Sigma \right)$ as inputs with $\Sigma$ a non-diagonal covariance matrix is a member of the class of functions defined by    
$$
\mathcal{D}_{d, n} =\left\{ \M : \R^d \to \R^\NN\, : \; \bo{\X} \sim F,\; F(\bo{x}) \neq \prod_{j=1}^d F_j(x_j);\; \bo{x} \in \R^d  \right\} \, .  
$$   
               
\subsection{New insight into dependency functions} \label{sec:mdep}  
In this section, we recall useful results about generic dependency models of non-independent variables (see \cite{skorohod76,lamboni22b,lamboni21,lamboni21ar,lamboni22mcap}). For a $d$-dimensional random vector of non-independent variables (i.e., $\bo{\X} \sim F$), a dependency model consists in expressing a subset of $d-1$ variables as a function of independent variables. \\   
       
Formally, if $\bo{\X} \sim F$ with $F(\bo{x}) \neq \prod_{j=1}^d F_j(x_j)$, then there exists  (\cite{skorohod76,lamboni22b,lamboni21,lamboni21ar,lamboni22mcap}) \\  
$\quad$ (i) new independent variables $\bo{Z} := \left(Z_{w_1}, \ldots, Z_{w_{d-1}}\right)$, which are independent of $\X_j$; \\
$\quad$ (ii) a dependency function $r_j  : \R^d  \to R^{d-1}$; \\ 
such that     
\begin{equation}    
\bo{\X}_{\sim j} \stackrel{d}{=} r_j\left(\X_j, \bo{Z}\right);
 \qquad \mbox{and} \qquad 
\left(\X_j, \,\bo{\X}_{\sim j}\right) \stackrel{d}{=} \left(\X_j, \, r_j\left(\X_j, \bo{Z}\right) \right) \, , 
\end{equation}    
where $\bo{\X}_{\sim j} =: (\X_{w_1}, \ldots, \X_{w_{d-1}})$, and $A\stackrel{d}{=}B$ means that the random variables $A$ and $B$ have the same CDF.\\ 
 
 It is worth noting that the dependency model is not unique in general. The uniqueness can be obtained under additional conditions given in  Proposition \ref{prop:uniq}, which enable the inversion of the dependency function $r_j$.     
        
\begin{prop}  \label{prop:uniq} 
Consider a dependency model of the continuous random vector $\bo{\X} \sim F$ given by $\bo{\X}_{\sim j} = r_j\left(\X_j, \bo{Z}\right)$ with a prescribed order $\left(w_1, \ldots, w_{d-1} \right)$. \\     
 
If $\X_j$ is the explanatory variable and the distribution of $\bo{Z}$ is prescribed, then  \\ 
  
$\quad$ (i) the dependency model is uniquely defined. \\ 
    
$\quad$ (ii) the dependency model is invertible and the unique inverse is given by
\begin{equation} \label{eq:invdep}
 \bo{Z} = r_{j}^{-1}\left(\bo{\X}_{\sim j}\,|\, \X_{j}\right) \, .
\end{equation}  
\end{prop}       
\begin{preuve} 
See Appendix \ref{app:prop:uniq}.   
\begin{flushright}  
$\Box$         
\end{flushright} 
\end{preuve}            
   
It is to be noted that the dependency models are vector-valued functions of independent input variables, which ease the calculus of partial derivatives, including partial derivatives of $\bo{\X}$ w.r.t. $\X_j$. Moreover, the inverse of the dependency model avoids working with $\bo{\Z}$. A natural choice of the order $\left(w_1, \ldots, w_{d-1} \right)$ is $(1, \ldots, j-1, j+1, \ldots, d)$. 
										        		 			  
\subsection{Enhance implicit functions: dependency functions} \label{sec:expl}
In this section, we provide a probabilistic version of the implicit functions using dependency models.  \\   

Consider $\bo{\X} \sim F$, a sample value or a realization of $\bo{\X}$ given by $\bo{x} \in \R^d$  and a function $h : \R^d \to \R^p$ with $p\leq d$. When connecting the input variables by $p$ compatible equations, that is, 
$$
h(\bo{\X}) = \bo{0} \, ,  
$$ 
the theorem of implicit function claims that for each sample value $\bo{x}^*$ satisfying $h(\bo{x}^*)=\bo{0}$, a subset of $\bo{\X}$ can be expressed as a function of the others in the neighborhood of $\bo{x}^*$. To recall such theorem, we use $u\subseteq \{1, \ldots, d\}$ with $|u| := card(u)=d-p$, $\bo{\X}_u :=\left(\X_j, \, \forall\,  j\in u \right)$, $\bo{\X}_{\sim u} := \left(\X_j, \, \forall\, j \in \{1,\ldots d\} \setminus u \right)$ and $B(\bo{x}^*_{u}, r_1) \subseteq \R^{d-p}$ (resp. $B(\bo{x}^*_{\sim u}, r_2)$) for an open ball centered on $\bo{x}^*_{u}$ (resp. $\bo{x}^*_{\sim u}$) with radius $r_1$ (resp. $r_2$). Again, $h_{\bo{x}_{\sim u}}'(\bo{x}^*)$  (resp. $h_{\bo{x}_{u}}'$) is the  formal Jacobian of $h$ w.r.t.. to $\bo{x}_{\sim u}$ (resp. $\bo{x}_{u}$).\\        
     
\begin{theorem} (implicit function) \label{theo:imp} \\
Assume that $h(\bo{x}^*)=\bo{0}$ and $h_{\bo{x}_{\sim u}}'(\bo{x}^*)$ is invertible. 
Then,  there exists a function $g : B(\bo{x}^*_{u}, r_1) \to B(\bo{x}^*_{\sim u}, r_2)$ such that 
$$
\bo{x}_{\sim u} = g(\bo{x}_{u}); \quad \qquad   \bo{x}_{\sim u}^* = g(\bo{x}_{u}^*);
\qquad \quad   
 \frac{\partial \bo{x}_{\sim u}}{\partial \bo{x}_{u}}(\bo{x}_{u}^*) = -\left[h_{\bo{x}_{\sim u}}'(\bo{x}^*)\right]^{-1} h_{\bo{x}_{u}}'(\bo{x}^*)  \, .  
$$  
\end{theorem}       
       
While Theorem \ref{theo:imp} is useful, it turns out that the implicit function theorem gives a local relationship among the variables. Remark that the dependency models derived in Section \ref{sec:mdep} provide the global relationships once the CDFs of the input variables are known. The distribution function of the variables that complies with the constraints given by $h(\bo{\X}) = \bo{0}$ is needed for building a global implicit function. To derive such distribution function, we assume that  
      
(A1) all the constraints $h(\bo{\X})=\bo{0}$ are compatible. \\   
   
Under (A1), the constraints $h(\bo{\X})=\bo{0}$ introduce new dependency structures on the initial CDF $F$, which mater for our analysis. The probability theory ensures the existence of a distribution function that captures such dependencies. 
    
\begin{prop} \label{prop:cdist}    
Let $\bo{\X} \sim F$ and $\bo{\X}^c :=\left\{\bo{\X} \sim F \, : \, h(\bo{\X}) = \bo{0} \right\} \sim F^c$ be the constrained variables. If (A1) holds, then we have
\begin{equation}  
\displaystyle 
\left\{ \begin{array}{l} 
 \bo{\X} \sim F  \\ 
 \mbox{s.t.} \quad h(\bo{\X}) := \bo{0} \\ 
\end{array} \right.  \stackrel{d}{=}  \bo{\X}^c\, ,      \nonumber   
\end{equation}     
where $\stackrel{d}{=}$ denotes the equality in distribution. 
\end{prop}       
           
It comes out that introducing constraints on initial variables leads us to work with constrained variables following a new CDF, that is, $\bo{\X}^c \sim F^c$. Some generic and constrained variables $\bo{\X}^c$ and their corresponding distribution functions can be found in \cite{lamboni21ar,durante22,lamboni22b}. When analytical derivations of the CDF of $\bo{\X}^c$ are hard or impossible, a common practice consists in fitting a distribution function to the observations of $\bo{\X}^c$ by means of numerical simulations (See \cite{rosenblatt56,parzen62,epanechnikov69,lamboni21ar,mcneil15,durante15} for examples of the distributions and  densities estimations). Using the new distributions of the input variables, Corollary \ref{coro:depexp} provides the probabilistic version of the implicit function.  
              
\begin{defi}         
A distribution $G$ is said to be a degenerate CDF whenever $G$ is the CDF of the Dirac measure having $\delta_{a}$ as the probability mass function with $a \in \R$. 
\end{defi}    
           
\begin{corollary}  \label{coro:depexp}  
Consider a random vector $\bo{\X}^c :=\left\{\bo{\X} \sim F \, : \, h(\bo{\X}) = \bo{0} \right\}$ that follows $F^c$ as CDF. Assume (A1) holds and $F^c$ is a non degenerate CDF.\\
   
Then,  there exists a function $r_j : \R^{d} \to \R^{d-1}$ and $d-1$ independent variables $\bo{Z} \sim F_{\bo{Z}}$ such that  $\bo{Z}$ is independent of $\X_{j}^c$ and 
$$
\bo{\X}_{\sim j}^c = r_j(\X_{j}^c, \, \bo{Z}) \, .   
$$       
\end{corollary} 
\begin{preuve}
Such result is the dependency model for the distribution $F^c$ (see Section \ref{sec:mdep}).   
\begin{flushright}  
$\Box$
\end{flushright}
\end{preuve}

While Corollary \ref{coro:depexp} gives the explicit function that links $\bo{\X}_{\sim j}$ to $\X_j$, we can sometime extend such result as follows:
\begin{equation}
\bo{\X}_{\sim u}^c = r_u(\bo{\X}_{u}^c, \, \bo{Z}_u) \, , 
\end{equation}     
where $\bo{Z}_u$ is a vector of $d-|u|$ independent variables, which are independent of $\bo{\X}_{u}^c$ (see Section \ref{sec:mdep} and \cite{lamboni21ar}). 
        
\begin{rem}
We can easily generalize the above process to cope with i) constrained inequations such as $h(\bo{\X}) <\bo{0}$ or $h(\bo{\X}) > \bo{0}$, ii) a mixture of constrained equations and inequations involving different variables. 
\end{rem}              

\subsection{Representation of function with non-independent variables} \label{sec:expl}
In general a function may include a group of independent variables as well as groups of non-independent variables such as correlated variables and/or dependent variables. We can then organize such input variables as follows:\\ 
   
(O):  the random vector $\bo{\X} := (\X_1, \ldots, \X_d)$ is consisted of $K$ independent random vector(s) given by $\bo{\X}=\left(\bo{\X}_{\boldsymbol{\pi}_1}, \ldots, \bo{\X}_{\boldsymbol{\pi}_K} \right)$ where the sets $\boldsymbol{\pi}_1, \ldots, \boldsymbol{\pi}_K$ form a partition of $\{1, \ldots, d\}$.  
The random vector $\bo{\X}_{\boldsymbol{\pi}_{k_1}} := \left(\X_\jmath, \, \forall \, \jmath \in \boldsymbol{\pi}_{k_1} \right)$ is independent of $\bo{\X}_{\boldsymbol{\pi}_{k_2}} :=\left(\X_\jmath, \, \forall \, \jmath \in \boldsymbol{\pi}_{k_2} \right)$ for every  pair $k_1, k_2 \in \{1, \ldots, K\}$ with $k_1 \neq k_2$. Without loss of generality, we use $\bo{\X}_{\boldsymbol{\pi}_1}$ for a random vector of $d_1\geq 0$ independent variable(s); $\bo{\X}_{\boldsymbol{\pi}_k}$ with $k \geq 2$ for a random vector of $d_k\geq 2$ dependent variables.\\  
      
We use $(\pi_{1,k}, \ldots, \pi_{d_k,k})$ for the ordered permutation of $\boldsymbol{\pi}_k$ (i.e., $\pi_{1,k} < \pi_{2,k} \ldots < \pi_{d_k,k}$). For any $\pi_{j,k} \in \boldsymbol{\pi}_k$, we use $\X_{\pi_{j,k}}$ for an element of $\bo{\X}_{\boldsymbol{\pi}_k}$; $(\sim\pi_{j,k}) := (\pi_{1,k}, \ldots, \pi_{j-1,k}, \pi_{j+1,k}, \ldots, \pi_{d_k,k})$ and $\bo{Z}_k :=\left(Z_{\pi_{\ell, k}}, \, \forall \, \pi_{\ell, k} \in (\sim \pi_{j,k})  \right)$.              
Using the dependency function (see Section \ref{sec:mdep}), we can represent $\bo{\X}_{\boldsymbol{\pi}_k}$ by     
\begin{equation} \label{eq:depmxk} 
    \bo{\X}_{\sim \pi_{j,k}} \stackrel{d}{=} r_{\pi_{j,k}}\left(\X_{\pi_{j,k}},\bo{Z}_k \right), \quad \forall\, k \in \{2,\ldots, K\}		\, ,          
\end{equation}    
where  $r_{\pi_{j,k}}(\cdot) =\left(r_{\pi_{1,k}}(\cdot), \ldots, r_{\pi_{j-1,k}}(\cdot), \; \; r_{\pi_{j+1,k}}(\cdot), \ldots, r_{\pi_{d_k}}(\cdot) \right)$;       
$\bo{Z}_k$ is a random vector of $d_k-1$ independent variables, and $\X_{\pi_{1,k}}$  is independent of $\bo{Z}_k$.     
Based on the above dependency model of $\bo{\X}_{\boldsymbol{\pi}_k}$ with $k=2, \ldots, K$, let us introduce new functions $c_{\pi_{j,k}} :\R^{d_k} \to \R^{d_k}$ given by 
$ 
c_{\pi_{j,k}}(\cdot) = \left(r_{\pi_{1,k}}(\cdot), \ldots, r_{\pi_{j-1,k}}(\cdot),\;\;  r_{\pi_{j,k}}(\cdot) =\X_{\pi_{j,k}},\; \; r_{\pi_{j+1,k}}(\cdot), \ldots, r_{\pi_{d_k}}(\cdot) \right)   
$
and     
$$
\bo{\X}_{\boldsymbol{\pi}_k} \stackrel{}{=} 
c_{\pi_{j,k}}(\X_{\pi_{j,k}},\, \bo{Z}_k) =: \left(\X_{\pi_{j,k}},\, r_{\pi_{j,k}}(\X_{\pi_{j,k}}, \, \bo{Z}_k)  \right) \, .  
$$     
The function $c_{\pi_{j,k}}$ maps independent variables $(\X_{\pi_{j,k}},\, \bo{Z}_k)$ onto $\bo{\X}_{\boldsymbol{\pi}_k}$, and the chart    
\begin{equation}     \label{eq:comp}        
\displaystyle   
\begin{array}{lcccccc} & \R^{d} &  \stackrel{c}{\to} & \R^d & \stackrel{\M}{\to} & \R^\NN \\
                    &\left(\begin{array}{c} \bo{\X}_{\boldsymbol{\pi}_1} \\
									  \X_{\pi_{j,2}} \\ \bo{Z}_{2} \\ \vdots \\
										\X_{\pi_{j,K}} \\ \bo{Z}_{K} 
										\end{array}\right)
										& \mapsto & \left(\begin{array}{c} \bo{\X}_{\boldsymbol{\pi}_1} \\ 
											\bo{\X}_{\boldsymbol{\pi}_2} \\ \vdots \\
	        		      \bo{\X}_{\boldsymbol{\pi}_K} \\ \end{array}\right) =\bo{\X}				
									& \mapsto & \left(\begin{array}{c} \M_1(\bo{\X})\\ 
									   \vdots \\ 	\M_\NN(\bo{\X}) \\
										\end{array}\right)
										\end{array}   \,  ,  \nonumber    
\end{equation} 
leads to a new representation of functions with non-independent variables by composing $\M$ by $c$, that is,  
\begin{equation} \label{eq:repfniid} 
\M(\bo{\X}_{\boldsymbol{\pi}_1}, \, \bo{\X}_2, \ldots, \bo{\X}_{\boldsymbol{\pi}_k}) \stackrel{d}{=} \M \circ c\left(\bo{\X}_{\boldsymbol{\pi}_1},\, \X_{\pi_{j,2}},\, \bo{Z}_{2}, \ldots, \X_{\pi_{j,K}},\, \bo{Z}_{K} \right) \, ,       
\end{equation}            
with $c\left(\bo{\X}_{\boldsymbol{\pi}_1},\, \X_{\pi_{j,2}},\, \bo{Z}_{2}, \ldots, \X_{\pi_{j,K}},\, \bo{Z}_{K} \right) := \left(\bo{\X}_{\boldsymbol{\pi}_1},\, c_{\pi_{j,2}}(\X_{\pi_{j,2}},\, \bo{Z}_2), \ldots, c_{\pi_{j,K}}(\X_{\pi_{j,K}},\, \bo{Z}_K) \right)$.\\
   
The equivalent representation of $\M$ given by (\ref{eq:repfniid}) relies on the innovation variables $\bo{Z} := (\bo{Z}_2, \ldots, \bo{Z}_K)$. Recall that for the continuous random vector $\bo{\X}_{\boldsymbol{\pi}_k}$, the dependency model $r_{\pi_{j,k}}$ given by (\ref{eq:depmxk}) is always invertible (see Proposition \ref{prop:uniq}), and therefore $c_{\pi_{j,k}}$ is also invertible. Such inversions are helpful for working with $\bo{\X}$ only.   
                   
\section{Actual partial derivatives} \label{sec:fwni} 
This section deals with the calculus of partial derivatives of function with non-independent variables using only one relationship among inputs such as the dependency model given by Equation (\ref{eq:depmxk}). The usual assumptions made are \\
  
(A2) the joint (resp. marginal) CDF is continuous and has a density function $\rho >0$ on its open support; \\  

(A3) each component of the dependency function $r_{\pi_{j,k}}$ is differentiable w.r.t. $\X_{\pi_{j,k}}$; \\ 
   
(A4) each component of the function $\M$, that is, $\M_\ell$ with $\ell=1,\,\ldots,\,\NN$ is differentiable w.r.t. each input.  \\         

Without loss of generality and for the sequel of simplicity, we suppose that $\NN=1$ in what follows. Namely, we use $\mathcal{I}_{d\times d} \in \R^{d \times d}$ for the identity matrix; $\mathsf{O}_{d\times d_1} \in \R^{d \times d_1}$ for a null matrix. 
It is common to use $\nabla\M_{\bo{x}_{\boldsymbol{\pi}_k}}$ for the formal partial derivatives of $\M$ w.r.t. each input of $\bo{\X}_{\boldsymbol{\pi}_k}$ (i.e., the derivatives obtained by considering inputs as independent) with $k=1, \ldots, K$. Thus, the formal gradient of $\M$ (i.e., the gradient w.r.t. the Euclidean metric) is given by  
\begin{equation} \label{eq:fgrdf}  
\nabla \M := \left[\nabla\M_{\bo{x}_{\boldsymbol{\pi}_1}}^\T \; \nabla\M_{\bo{x}_{\boldsymbol{\pi}_2}}^\T  \ldots \nabla\M_{\bo{x}_{\boldsymbol{\pi}_K}}^\T \right]^T \, . 
\end{equation}     
        
Keeping in mind the function $c_{\pi_{j,k}}(\cdot)$, the partial derivatives of each component of $\bo{\X}_{\boldsymbol{\pi}_k}$ w.r.t. $\X_{\pi_{j,k}}$ are given by   
\begin{equation} \label{eq:fgrdc}      
J^{(\pi_{j,k})} :=  \frac{\partial c_{\pi_{j,k}}}{\partial x_{\pi_{j,k}}} =  \left[\frac{\partial \X_{\pi_{1,k}}}{\partial x_{\pi_{j,k}}}\, \ldots \, \frac{\partial \X_{\pi_{d_k,k}}}{\partial x_{\pi_{j,k}}} \right]^T = \left[\frac{\partial r_{\pi_{1,k}}}{\partial x_{\pi_{j,k}}}\, \ldots \, \underbrace{1}_{j^{\text{th}} \, \text{position}} \, \ldots \, \frac{\partial r_{\pi_{d_k,k}}}{\partial x_{\pi_{j,k}}} \right]^T\,  .   
\end{equation}            
We use $J^{(\pi_{j,k})} _i$  for the $i^{\mbox{th}}$ element of $J^{(\pi_{j,k})}$. For instance, $J^{(\pi_{j,k})}_{j} =1$ and $J^{(\pi_{j,k})}_{d_k}$ represents the partial derivative of $\X_{\pi_{d_k,k}}$ w.r.t. $\X_{\pi_{j,k}}$.  
It is worth recalling that $J^{(\pi_{j,k})}$ is a vector-valued function of $(\X_{\pi_{j,k}}, \bo{Z}_k)$, and Lemma \ref{lem:unigrad} expresses $J^{(\pi_{j,k})}$ as a function of $\bo{x}_{\boldsymbol{\pi}_k}$ only. 

\begin{lemma} \label{lem:unigrad} 
Let  $\bo{x}_{\boldsymbol{\pi}_k}$ be a sample value of $\bo{\X}_{\boldsymbol{\pi}_k}$. If assumptions (A2)-(A4) hold, then the partial derivatives of $\bo{\X}_{\boldsymbol{\pi}_k}$ w.r.t. $\X_{\pi_{j,k}}$ evaluated at $\bo{x}_{\boldsymbol{\pi}_k}$ is given by  
\begin{equation} \label{eq:unigrad}    
J^{(\pi_{j,k})}\left(\bo{x}_{\boldsymbol{\pi}_k} \right) := \left[\frac{\partial r_{\pi_{1,k}}}{\partial x_{\pi_{j,k}}}\, \ldots \, \underbrace{1}_{j^{\text{th}} \, \text{position}} \, \ldots \, \frac{\partial r_{\pi_{d_k,k}}}{\partial x_{\pi_{j,k}}} \right]^T\left(x_{\pi_{j,k}},\, r_{\pi_{j,k}}^{-1}\left(\bo{x}_{\sim \pi_{j,k}} |x_{\pi_{j,k}} \right) \right)\,  .    
\end{equation}                     
\end{lemma}                   
\begin{preuve}     
See Appendix \ref{app:lem:unigrad}. 
\begin{flushright}     
$\Box$
\end{flushright}
\end{preuve} 

Again, $J^{(\pi_{j,k})}_\ell \left(\bo{x}_{\boldsymbol{\pi}_k} \right) $ with $\ell =1, \ldots, d_k$ stands for the $\ell^{\mbox{th}}$ component  of $J^{(\pi_{j,k})}\left(\bo{x}_{\boldsymbol{\pi}_k} \right)$ provided in Lemma \ref{lem:unigrad}. Using such components and the chain rule,   
Theorem \ref{theo:agrad} provides the actual partial derivatives of functions with non-independent variables (i.e., $\frac{\partial_a \M}{\partial \bo{x}}$), that is, the derivatives obtained by making use of only one dependency function given by Equation (\ref{eq:depmxk}).  
          
\begin{theorem} \label{theo:agrad}  Let $\bo{x} \in \R^d$ be a sample value of $\bo{\X}$ and $\pi_{j,k} \in \boldsymbol{\pi}_k$ with $k=2, \ldots, K$.  
If assumptions (A2)-(A4) hold, then \\  
  
$\quad$ (i) the actual Jacobian matrix of $c_{\pi_{j,k}}$ or $\frac{\partial_a \bo{\X}_{\boldsymbol{\pi}_k}}{\partial \bo{x}_{\boldsymbol{\pi}_k}}$ is given by  
\begin{equation} \label{eq:ajacxk}  
J_{\boldsymbol{\pi}_k}^a \left(\bo{x}_{\boldsymbol{\pi}_k}\right) := \left[
\begin{array}{cccc}
\frac{J^{(\pi_{j,k})} \left(\bo{x}_{\boldsymbol{\pi}_k} \right)}{J^{(\pi_{j,k})}_1\left(\bo{x}_{\boldsymbol{\pi}_k} \right)} & \ldots & \frac{J^{(\pi_{j,k})}\left(\bo{x}_{\boldsymbol{\pi}_k} \right)}{J^{(\pi_{j,k})}_{d_k}\left(\bo{x}_{\boldsymbol{\pi}_k} \right)} \\  
\end{array}   
\right], \quad \forall\, k \in \{2, \ldots, K\}  \, .  
\end{equation}   
            
$\quad$ (ii) The actual Jacobian matrix of $c$ or $\frac{\partial_a \bo{\X}}{\partial \bo{x}}$ is given by   
\begin{eqnarray}    \label{eq:ajacall} 
J^a(\bo{x})  &:=& \left[      
\begin{array}{cccc}
\mathcal{I}_{d_1\times d_1} & \mathsf{O}_{d_1 \times d_2} & \ldots & \mathsf{O}_{d_1 \times d_K} \\ 
\mathsf{O}_{d_2 \times d_1} & J_{\boldsymbol{\pi}_2}^a\left(\bo{x}_{\boldsymbol{\pi}_2}\right) & \ldots &  \mathsf{O}_{d_2 \times d_K} \\
\vdots & \vdots  & \ddots & \vdots \\  
\mathsf{O}_{d_K \times d_1} & \mathsf{O}_{d_K \times d_2} & \ldots & J_{\boldsymbol{\pi}_K}^a\left(\bo{x}_{\boldsymbol{\pi}_K} \right)\\
\end{array}  
 \right]   \, ,               
\end{eqnarray}       

$\quad$ (iii) The actual partial derivatives of $\M$ are given by 
\begin{equation}  \label{eq:agradf}
\frac{\partial_a \M}{\partial \bo{x}}(\bo{x})  := J^a \left(\bo{x} \right)^\T \nabla \M(\bo{x})  \, .   
\end{equation}               
\end{theorem}   
\begin{preuve}  
see Appendix \ref{app:theo:agrad}.  
\begin{flushright}      
$\Box$  
\end{flushright}     
\end{preuve}   
  
Results from Theorem \ref{theo:agrad} are based only on one dependency function, which uses $\X_{\pi_{j,k}}$ as the explanatory input. Thus, the actual Jacobian $J^a(\bo{x})$ and the actual partial derivatives of $\M$ provided in (\ref{eq:ajacxk})-(\ref{eq:agradf}) are going to change with the choice of the explanatory input $\X_{\pi_{j,k}}$ for every $j \in \{1,\ldots, d_k\}$. All these possibilities are not surprising. Indeed, while no additional explicit assumption is necessary for calculating the partial derivatives of $\bo{\X}_{\boldsymbol{\pi}_k}$ w.r.t. $\X_{\pi_{j,k}}$ (i.e., $J^{(\pi_{j,k})}$), we implicitly keep the other variables fixed when calculating the partial derivative of $\bo{\X}_{\boldsymbol{\pi}_k}$ w.r.t. $\X_{\pi_{i,k}}$, that is, $\frac{J^{(\pi_{j,k})}}{J^{(\pi_{j,k})}_i}$ for each $i\neq j$. Such implicit assumption is due to the reciprocal rule used to derive the results (see Appendix \ref{app:theo:agrad}). In general, the components of $\bo{\X}_{\sim \pi_{j,k}}$ such as $\X_{\pi_{i_1,k}}$ and $\X_{i_2,k}$ are both function of $\X_{\pi_{j,k}}$ and $Z_{\pi_{1,k}}$ at least. Thus, different possibilities of actual Jacobians are based on different implicit assumptions, and it becomes difficult to use the actual partial derivatives. Further drawbacks of the actual partial derivatives of $\M$ are illustrated in Example 1. 
       
\subsection*{\textbf{Example 1}}
We consider the function  
$  
\M(\X_1,\, \X_2) = \X_1 +\X_2 +  \X_1 \X_2 
$, 
which includes two correlated inputs $\bo{\X} \sim \mathcal{N}_2\left(\bo{0},\, \left[\begin{array}{cc}1 &\rho\\ \rho & 1  \end{array}\right]\right)$. We see that 
$
\nabla \M(\bo{X}) =  \left[\begin{array}{c} 1+ X_2\\1+ \X_1  \end{array}\right] 
$.           
Using the dependency model of $\bo{X}$ given by  (see \cite{lamboni21,lamboni21ar,lamboni22b,lamboni22d})         
$$          
X_2=\rho X_1 + \sqrt{1-\rho^2} Z_2  \Longleftrightarrow Z_2 = (X_2 - \rho X_1)/	\sqrt{1-\rho^2} \, ,
$$   
 the actual Jacobian matrix of $c$ and the actual partial derivatives of $\M$ are given by 
$$    
J^a(\bo{\X}) = \left[\begin{array}{cc}1 & \frac{1}{\rho} \\ \rho & 1  \end{array}\right];  
\qquad
\frac{\partial_a \M}{\partial \bo{x}}(\bo{\X}) = \left[\begin{array}{c}1 +\X_2+\rho(1+\X_1) \\ \frac{1 +\X_2}{\rho} + 1+\X_1 \end{array}\right] \, .  
$$ 
When $\rho = 1$, both inputs are perfectly correlated, and we have $\X_1 =\X_2$, which also implies $\M(\X_1, \X_2) = \M(\X_1) =2\X_1 + \X^2_1 = \M(\X_2) =2\X_2 + \X^2_2$. We can check that $\frac{\partial_a \M}{\partial \bo{x}}(\bo{\X}) = \left[\begin{array}{c}2 +2 \X_1\\ 2 +2 \X_2 \end{array}\right]$. However,  When $\rho = 0$, both inputs are independent, and we should expect the actual partial derivatives to be equal to the formal gradient $\nabla \M$, but it is not the case.  
Moreover, using the second dependency model, that is,           
$$          
X_1=\rho X_2 + \sqrt{1-\rho^2} Z_1  \Longleftrightarrow Z_1 = (X_1 - \rho X_2)/	\sqrt{1-\rho^2} \, ,  
$$ 
it comes out that 
$   
J^a(\bo{\X}) = \left[\begin{array}{cc}1 & \rho \\ \frac{1}{\rho} & 1  \end{array}\right];
\qquad 
\frac{\partial_a \M}{\partial \bo{\X}}(\bo{\X}) = \left[\begin{array}{c}1 +\X_2+\frac{1+\X_1}{\rho} \\ \rho(1 +\X_2) + 1+\X_1 \end{array}\right] \, ,     
$            
which differ from the previous results. All these drawbacks are due to the implicit assumptions  made, which can be avoided (see Section \ref{sec:dgrad}). 
                  
\section{Dependent Jacobian and partial derivatives} \label{sec:dgrad}
 This section aims at deriving the  first and second-order partial derivatives of functions with non-independent variables without using any additional assumption neither explicit nor implicit. 
Basically, we are going to calculate or compute the partial derivatives of $\bo{\X}_{\boldsymbol{\pi}_k}$ w.r.t. $\X_{\pi_{i,k}}$ using only the dependency function that includes $\X_{\pi_{i,k}}$ as explanatory input, that is, 
$$
\bo{\X}_{\sim \pi_{i,k}} = r_{\pi_{i,k}}\left(\X_{\pi_{i,k}}, \bo{\Z}_{\pi_{i,k}} \right) \, , \quad
\forall \, i =1, \ldots, d_k;   
\quad  
\forall \;  k = 2, \ldots, K \, .        
$$       
Using the above dependency function, the partial derivatives of $\bo{\X}_{\boldsymbol{\pi}_k}$ w.r.t. $\X_{\pi_{i,k}}$ is given as follows (see (\ref{eq:unigrad})):
\begin{equation} \label{eq:grdxk}    
J^{(\pi_{i,k})}\left(\bo{x}_{\boldsymbol{\pi}_k} \right) :=
\frac{\partial \bo{\X}_{\boldsymbol{\pi}_k}}{\partial x_{\pi_{i,k}}} = 
 \left[\frac{\partial r_{\pi_{1,k}}}{\partial x_{\pi_{i,k}}}\, \ldots \, \underbrace{1}_{i^{\text{th}} \, \text{position}} \, \ldots \, \frac{\partial r_{\pi_{d_k,k}}}{\partial x_{\pi_{i,k}}} \right]^T\left(x_{\pi_{i,k}},\, r_{\pi_{i,k}}^{-1}\left(\bo{x}_{\sim \pi_{i,k}} |x_{\pi_{i,k}} \right) \right)\,  .    \nonumber         
\end{equation}              
 It is to be noted that $J^{(\pi_{i,k})}$ does not require any supplementary assumption, as $\X_{\pi_{i,k}}$ and $\bo{\Z}_{\pi_{i,k}}$ are independent.
Thus, $d_k$ different dependency models are necessary to derive the dependent Jacobian  and partial derivatives of $\M$ (see Theorem \ref{theo:dgrad}). 
     
\begin{theorem} \label{theo:dgrad}  Let $\bo{x} \in \R^d$ be a sample value of $\bo{\X}$, and assume that (A2)-(A4) hold.\\        
   
$\quad$ (i) For all $k\geq 2$, the dependent Jacobian matrix $\frac{\partial \bo{\X}_{\boldsymbol{\pi}_k}}{\partial \bo{x}_{\boldsymbol{\pi}_k}}$ is given by 
\begin{equation} \label{eq:djacxk}  
J_{\boldsymbol{\pi}_k}^d \left(\bo{x}_{\boldsymbol{\pi}_k} \right) := \left[ 
\begin{array}{cccc}
J^{(\pi_{1,k})} \left(\bo{x}_{\boldsymbol{\pi}_k} \right) & \ldots & J^{(\pi_{d_k,k})} \left(\bo{x}_{\boldsymbol{\pi}_k} \right) \\     
\end{array}       
\right]  \, .        
\end{equation}           
       
$\quad$ (ii) The dependent Jacobian matrix $\frac{\partial \bo{\X}}{\partial \bo{x}}$  is given by   
\begin{eqnarray}    \label{eq:djacall} 
J^d(\bo{x})  &:=& \left[      
\begin{array}{cccc}
\mathcal{I}_{d_1\times d_1} & \mathsf{O}_{d_1 \times d_2} & \ldots & \mathsf{O}_{d_1 \times d_K} \\ 
\mathsf{O}_{d_2 \times d_1} & J_{\boldsymbol{\pi}_2}^d\left(\bo{x}_{\boldsymbol{\pi}_2}\right) & \ldots &  \mathsf{O}_{d_2 \times d_K} \\
\vdots & \vdots  & \ddots & \vdots \\  
\mathsf{O}_{d_K \times d_1} & \mathsf{O}_{d_K \times d_2} & \ldots & J_{\boldsymbol{\pi}_K}^d\left(\bo{x}_{\boldsymbol{\pi}_K} \right)\\
\end{array}   
 \right]   \, ,           
\end{eqnarray}       
 
$\quad$ (iii) The partial derivatives of $\M$ are given by 
\begin{equation}  \label{eq:dgradf}
\frac{\partial \M}{\partial \bo{x}}(\bo{x})  :=  J^d \left(\bo{x} \right)^\T  \nabla \M(\bo{x})   \, , 
\end{equation} 
\end{theorem}       
\begin{preuve}  
See Appendix \ref{app:theo:dgrad}. 
\begin{flushright} 
$\Box$ 
\end{flushright}
\end{preuve}   
 
Although the results from Theorem \ref{theo:dgrad} require different dependency models, such results are more comprehensive than the actual partial derivatives because
no supplementary assumption is available for each non-independent variables. \\  
      
To derive the second-order partial derivatives of $\M$, we use $\M_{x_{i} x_{j}}^{''}$ for the formal cross-partial derivative of $\M$ w.r.t. $x_{i},\, x_{j}$ and $H_{\boldsymbol{\pi}_k} := \left(\M_{x_{j_1} x_{j_2}}^{''} \right)_{j_1, j_2 \in \boldsymbol{\pi}_k}$  for the formal or ordinary Hessian matrices of $\M$ restricted to $\bo{\X}_{\boldsymbol{\pi}_k}$ with $k=1, \ldots, K$. In the same sense, we use $H_{\boldsymbol{\pi}_{k_1},\boldsymbol{\pi}_{k_2}} := \left(\M_{x_{j_{1}}  x_{j_2}}^{''} \right)_{j_1 \in \boldsymbol{\pi}_{k_1},\, j_2 \in \boldsymbol{\pi}_{k_2}}$ for the formal cross-Hessian matrix of $\M$ restricted to $(\bo{\X}_{\boldsymbol{\pi}_{k_1}}, \, \bo{\X}_{\boldsymbol{\pi}_{k_2}})$ for every pair $k_1, \, k_2 \in \{1, \ldots K\}$ with $k_1\neq k_2$.                       
To ensure the existence of the second-order partial derivatives, we assume that \\

(A5):  the function $\M$ is twice (formal) differentiable w.r.t. each input\\

(A6):  every dependency function $r_{\pi_{j,k}}$ is twice differentiable w.r.t. $\X_{\pi_{j,k}}$. \\   
   
By considering the $d_k$ dependency models of $\bo{\X}_{\boldsymbol{\pi}_k}$ (i.e., $\bo{\X}_{\sim \pi_{i,k}} =r_{\pi_{i,k}}\left(\X_{\pi_{i,k}},\, \bo{Z}_{\pi_{i,k}} \right)$ with $i=1, \ldots, d_k$) used to derive the dependent Jacobian, we can write
$$   
\frac{\partial J^{(\pi_{i,k})}}{\partial x_{\pi_{i,k}}}( \bo{x}_{\boldsymbol{\pi}_k}) := \frac{\partial^2 \bo{\X}_{\boldsymbol{\pi}_k} }{\partial^2 x_{\pi_{i,k}}} = \left[\frac{\partial^2 r_{\pi_{1,k}}}{\partial^2 x_{\pi_{i,k}}} \ldots \underbrace{0}_{i^{\text{th}} \, \text{position}} \; \ldots\; \frac{\partial^2 r_{\pi_{d_k,k}}}{\partial^2 x_{\pi_{i,k}}} \right]^T \left(x_{\pi_{i,k}},\, r_{\pi_{i,k}}^{-1}\left(\bo{x}_{\sim \pi_{i,k}} |x_{\pi_{i,k}} \right) \right) \, ,         
$$            
for the second partial derivatives of $\bo{\X}_{\boldsymbol{\pi}_k}$ w.r.t. $\X_{\pi_{i,k}}$. We then use  
\begin{equation} \label{eq;diffJhes}
\mathcal{D}J^d_{\boldsymbol{\pi}_k}\left(\bo{x} \right) := diag\left( \left[\frac{\partial J^{(\pi_{1,k})}}{\partial x_{\pi_{1,k}}}\left(\bo{x}_{\boldsymbol{\pi}_k} \right), \ldots,  \frac{\partial J^{(\pi_{d_k,k})}}{\partial x_{\pi_{d_k,k}}}\left(\bo{x}_{\boldsymbol{\pi}_k} \right) \right]^\T\nabla \M_{\bo{x}_{\boldsymbol{\pi}_k}}(\bo{x}) \right) J_{\boldsymbol{\pi}_k}^d\left(\bo{x}_{\boldsymbol{\pi}_k} \right) \, ,   
\end{equation}
for all $k \in \{2, \ldots, K\}$, and Theorem \ref{theo:dhess} provides the dependent second-order partial derivatives (i.e., $\frac{\partial^2 \M}{\partial^2 \bo{x}}$).        
  
   
\begin{theorem} \label{theo:dhess}  
Let $\bo{x}$ be a sample value of $\bo{\X}$.  If (A2), (A5) and (A6) hold, then        

$$
\frac{\partial^2 \M}{\partial^2 \bo{x}} (\bo{x}) := 
$$
\begin{eqnarray} 
\left[         
\begin{array}{cccc} 
H_{\boldsymbol{\pi}_1}(\bo{x}) & H_{\boldsymbol{\pi}_1, \boldsymbol{\pi}_2}(\bo{x}) J_{\boldsymbol{\pi}_2}^d\left(\bo{x}_{\boldsymbol{\pi}_2}\right) & \ldots & H_{\boldsymbol{\pi}_1, \boldsymbol{\pi}_K}(\bo{x}) J_{\boldsymbol{\pi}_K}^d\left(\bo{x}_{\boldsymbol{\pi}_k}\right)\\ 
 &  &  &  \\       
J_{\boldsymbol{\pi}_2}^d\left(\bo{x}_{\boldsymbol{\pi}_2}\right)^\T H_{\boldsymbol{\pi}_2, \boldsymbol{\pi}_1}(\bo{x}) & J_{\boldsymbol{\pi}_2}^d\left(\bo{x}_{\boldsymbol{\pi}_2}\right)^\T H_{\boldsymbol{\pi}_2}(\bo{x}) J_{\boldsymbol{\pi}_2}^d\left(\bo{x}_{\boldsymbol{\pi}_2}\right) & \ldots & J_{\boldsymbol{\pi}_2}^d\left(\bo{x}_{\boldsymbol{\pi}_2}\right)^\T H_{\boldsymbol{\pi}_2, \boldsymbol{\pi}_K}(\bo{x})  \\    
 & +\, \mathcal{D}J^d_{\boldsymbol{\pi}_2}\left(\bo{x} \right) &  & \times J_{\boldsymbol{\pi}_K}^d\left(\bo{x}_{\boldsymbol{\pi}_k}\right) \\    
\vdots & \vdots  & \ddots & \vdots \\    
J_{\boldsymbol{\pi}_K}^d\left(\bo{x}_{\boldsymbol{\pi}_k}\right)^\T H_{\boldsymbol{\pi}_K,\boldsymbol{\pi}_1}(\bo{x}) & J_{\boldsymbol{\pi}_K}^d\left(\bo{x}_{\boldsymbol{\pi}_K}\right)^\T H_{\boldsymbol{\pi}_K, \boldsymbol{\pi}_2}(\bo{x}) J_{\boldsymbol{\pi}_2}^d\left(\bo{x}_{\boldsymbol{\pi}_2}\right) & \ldots & J_{\boldsymbol{\pi}_K}^d\left(\bo{x}_{\boldsymbol{\pi}_K}\right)^\T H_{\boldsymbol{\pi}_K}(\bo{x}) \\     
&  &  & \times J_{\boldsymbol{\pi}_K}^d\left(\bo{x}_{\boldsymbol{\pi}_K}\right) + \,\mathcal{D}J^d_{\boldsymbol{\pi}_K}\left(\bo{x} \right) \\   
\end{array} \right]    \, .   \nonumber 
\end{eqnarray}     
\end{theorem}     
\begin{preuve}  
See Appendix \ref{app:theo:dhess}. 
\begin{flushright}    
$\Box$
\end{flushright}
\end{preuve}    
      
\subsection*{\textbf{Example 1 (revisited)}}
\noindent     
Since $                 
\M(\X_1,\, \X_2) = \X_1 +\X_2 +  \X_1 \X_2$ and the dependency models of $\bo{X}$ are given by 
$$          
X_2=\rho X_1 + \sqrt{1-\rho^2} Z_2  \Longrightarrow Z_2 = (X_2 - \rho X_1)/	\sqrt{1-\rho^2} \, , 
$$   
$$          
X_1=\rho X_2 + \sqrt{1-\rho^2} Z_1  \Longrightarrow Z_1 = (X_1 - \rho X_2)/	\sqrt{1-\rho^2} \, , 
$$ 
we can check that 
$$   
J^d = \left[\begin{array}{cc}1 & \rho \\ \rho & 1  \end{array}\right];
\qquad
\frac{\partial \M}{\partial \bo{x}}(\bo{\X}) = \left[\begin{array}{c}1 +\X_2+\rho(1+\X_1) \\ \rho(1 +\X_2) + 1+\X_1 \end{array}\right];
\quad 
\frac{\partial^2 \M}{\partial^2 \bo{x}}(\bo{X}) = \left[\begin{array}{cc} 2\rho & 1+\rho^2 \\ 1+\rho^2 & 2\rho \end{array}\right]   \, .
$$           
For instance, when $\rho = 1$, we have $\frac{\partial \M}{\partial \bo{x}}(\bo{\X}) = \frac{\partial_a \M}{\partial \bo{x}}(\bo{\X})$, and when $\rho=0$ we have $\frac{\partial \M}{\partial \bo{\X}}(\bo{X}) =\nabla \M(\bo{X})$ and  $\frac{\partial^2 \M}{\partial^2 \bo{\X}}(\bo{X}) =H(\bo{X})$. Thus, the dependent partial derivatives of $\M$ are coherent with the formal gradient and Hessian matrix when the inputs are independent. 
		  
\section{Expansion of functions with non-independent variables} \label{sec:ugrad} 
Although Section \ref{sec:dgrad} provides the partial derivatives and cross-partial derivatives of $\M$, it is misleading to think that infinitesimal increment of $\M$ given by $\M(\bo{\X}_{\boldsymbol{\pi}_k} + \epsilon \bo{e}_{\pi_{j,k}}) - \M(\bo{\X}_{\boldsymbol{\pi}_k})$ should result in an individual effect quantified by $\frac{\partial \M(\bo{\X}_{\boldsymbol{\pi}_k})}{\partial x_{\pi_{j,k}}} \epsilon$ with $\bo{e}_{\pi_{j,k}} := \left[0, \ldots,0, \underbrace{1}_{\pi_{j,k}^{\text{th}} \, \text{position}}, 0, \ldots, 0\right]^\T$ and $\epsilon >0$.  
Indeed,  moving $\X_{\pi_{j,k}}$ turns out to move the other variables partially, and the effects we obtain (i.e., $\nabla\M^\T_{\bo{x}_{\boldsymbol{\pi}_k}} J^{(\pi_{j,k})}$) can be imputed to other variables as well. The dependence structures of such effects are given by the dependent Jocobian matrix $J^{d}_{\boldsymbol{\pi}_k}(\bo{x}_{\boldsymbol{\pi}_k})$ (see Equation (\ref{eq:djacxk})). Therefore, the definition of gradient and the Hessian of $\M$ with non-independent variables need an introduction of the tensor metric or the Riemannian tensor.\\  
   
In differential geometry, the function of the form 
$$   
\begin{array}{cccl} 
c_{\pi_{j,k}} : & \R^{d_k}  & \to & \R^{d_k} \\
          & \left(\X_{\pi_{j,k}},\bo{Z}_{\pi_{j,k}} \right) & \mapsto & \bo{\X} := (\X_{\pi_{1,k}}, \ldots, \X_{\pi_{j,k}}, \ldots, \X_{\pi_{d_k,k}});
					\quad 
\bo{\X}_{\sim \pi_{j,k}} \stackrel{}{=} r_{\pi_{j,k}}\left(\X_{\pi_{j,k}}, \bo{Z}_{\pi_{j,k}} \right) \\  
\end{array}  \, ,      
$$      
for every $\pi_{i,k} \in\boldsymbol{\pi}_k$  can be seen as a parametrization of a manifold $\mathcal{M}_k$ in $\R^{d_k}$. The $d_k$ column entries of the dependent Jacobian matrix $J^{d}_{\boldsymbol{\pi}_k}(\bo{x}_{\boldsymbol{\pi}_k}) \in \R^{d_k \times d_k}$ span a local $m_k$-dimensional vector space a.k.a. the tangent space at $\bo{x}_{\boldsymbol{\pi}_k}$ with $m_k$ the rank of $J^{d}_{\boldsymbol{\pi}_k}(\bo{x}_{\boldsymbol{\pi}_k})$, that is,  the number of linearly independent columns of $J^{d}_{\boldsymbol{\pi}_k}(\bo{x}_{\boldsymbol{\pi}_k})$. \\ 
     
By considering all the $K$ groups of inputs and the corresponding dependent Jacobian matrix $J^d(\bo{x})$, we see that the support of the random vector $\bo{\X}$ is a $m$-dimensional manifold $\mathcal{M}$ in $\R^{d}$ with $m$ the rank of $J^d(\bo{x})$. When $m\leq d$, we are going to work in the tangent space $T\R^m$ (or local coordinate system) spanned out by the $m$ column entries of $J^d(\bo{x})$ that are linearly independent. Working in $T\R^m$ rather than $T\R^d$ ensures that the Riemannian tensor induced by $\bo{x} \mapsto \bo{x}$ using the dot product is invertible. Since the Riemannian tensor metric is often symmetric, the Moore-Penrose generalized inverse of symmetric matrices (\cite{moore1920,moore1935,penrose1955}) allows us to keep working in $T\R^d$ in what follows. Using the first fundamental form (see e.g., \cite{jost11,petersen16,sommer20}), the induced tensor metric is defined as   
\begin{equation} \label{eq:metric}       
G(\bo{x}) := J^d(\bo{x})^{\T}    J^d (\bo{x}) \, .     
\end{equation}       
         
Based on these elements, the gradient and the Hessian matrix are provided in Corollary \ref{coro:gradhess}. To that end, we use $G^{-1}$ for the inverse of the metric $G$ given by Equation (\ref{eq:metric}) when $m=d$ or the generalized  inverse of $G$ for every $m<d$ (\cite{moore1920,moore1935,penrose1955}). For any $k \in \{1, \ldots, d\}$, the Christoffel symbols are defined by (\cite{hanno82,jost11,sommer20,csaba22})  
$$
\Gamma^{k}_{i j} := \frac{1}{2} \sum_{\ell =1}^{m=d} G^{-1}_{k \ell}(\bo{x}) \left( 
 G_{i \ell, x_{j}}^{'} (\bo{x}) + 
G_{j \ell, x_{i}}^{'} (\bo{x}) - 
G_{i j, x_{\ell}}^{'} (\bo{x}) \right);\quad \forall\, i, j =1,\ldots, d \, ,
$$              
with $G_{i \ell, x_{j}}^{'}$ the formal partial derivative of $G_{i \ell}$ w.r.t. $ x_j$. 
          
\begin{corollary} \label{coro:gradhess} Let $\bo{x}$ be a sample value of $\bo{\X}$, and assume that (A2), (A5)-(A6) hold. \\    
   
$\quad$ (i) The gradient of $\M$ is given by   
\begin{equation}    \label{eq:grad}       
grad(\M)(\bo{x})  := G^{-1}(\bo{x}) \nabla \M(\bo{x})  \, ,         
\end{equation} 
         
$\quad$ (ii) The Hessian matrix of $\M$ is given by  
\begin{equation}  \label{eq:ugradf}
Hess_{i j}(\M)(\bo{x})  := \M_{x_{i} x_{j}}^{''}(\bo{x}) - 
\sum_{k=1}^{m=d}  \Gamma^{k}_{i j} (\bo{x})  \M_{x_{k}}^{'} (\bo{x})   \, . 
\end{equation}  
\end{corollary}          
\begin{preuve}    
Points (i)-(ii) result from the definition of the gradient and the Hessian matrix in a Riemannian geometric endowed with the metric $G$ (see \cite{jost11,sommer20,petersen16,vihua23}).
\begin{flushright} 
$\Box$     
\end{flushright}    
\end{preuve}    

Taylor expansion is widely used for approximating functions with independent variables. In what follows, we are concern with the approximation of a function with non-independent variables. The Taylor-type expansion of a function with non-independent variables is provided in Corollary \ref{coro:fexp} using the gradient and the Hessian matrix. 
   
\begin{corollary} \label{coro:fexp} Let $\bo{x}$, $\bo{x}_0$ be two sample values of $\bo{\X}$, and assume (A2), (A5)-(A6) hold. Then, we have
 
\begin{equation}    \label{eq:expf1}   
\M(\bo{x})   \approx  \M(\bo{x}_0)  +  \left(\bo{x}-\bo{x}_0\right)^\T grad(\M)(\bo{x}_0) + \frac{1}{2} \left(\bo{x}-\bo{x}_0\right)^\T Hess (\M)(\bo{x}_0) \left(\bo{x}-\bo{x}_0\right)  \, ,
\end{equation}     
provided that $\bo{x}$ is close to $\bo{x}_0$.          
\end{corollary}    
\begin{preuve}    
The proof is straightforward using the dot product induced by the tensor metric $G$ in the tangent space and bearing in mind the Taylor expansion provided in \cite{sommer20}. 
\begin{flushright}        
$\Box$       
\end{flushright}   
\end{preuve}         
         
\subsection*{\textbf{Example 1 (revisited)}}
For the function in Example 1, we can check that the tensor metric is \\ 
$ 
G(\bo{\X}) = \left[\begin{array}{cc}1 + \rho^2 &  2\rho \\ 2\rho & 1+ \rho^2   \end{array}\right];  
\qquad
G^{-1}(\bo{\X}) =\frac{1}{(\rho^2 -1)^2} \left[\begin{array}{cc}1 + \rho^2 &  -2\rho \\ -2\rho & 1+ \rho^2   \end{array}\right]
$; and  the gradient is   
$$
grad(\M)(\bo{\X})  = \frac{1}{(\rho^2 -1)^2}
\left[\begin{array}{c} 
 (1 - \rho)^2 + \X_2(1 + \rho^2) -2\rho \X_1 \\   
 (1 - \rho)^2 + \X_1(1 + \rho^2) -2\rho \X_2 
\end{array}\right] \, , 
$$             
which comes down to $\nabla\M(\bo{\X})$ when the variables are independent (i.e., $\rho=0$).  
          
\section{Conclusion} \label{sec:con}         
A new approach for calculating and computing the partial derivatives, gradient and Hessian of functions with non-independent variables is proposed and studied in this paper. It relies on i) dependency functions that model the dependency structures among dependent variables, including correlated variables, and ii) the tensor metric defined as the inner-product between the column entries of the Jocobian matrix of the dependency functions. Since the so-called dependent partial derivatives do not required any additional assumption (which is always the case), such derivatives should be preferred. The results obtained depend on the parameters of the distribution function or the density function of non-independent variables. For the values of such parameters that lead to independent variables, the proposed gradient comes down to the formal gradient or the gradient w.r.t. the Euclidean metric. \\ 
Using the proposed gradient and Hessian matrix, the Taylor-type expansion of a function with non-independent variables is provided. Although, the generalized inverse of a symmetric matrix is used in this paper, more investigation of the gradient calculus is needed when the tensor metric is not invertible. The proposed gradient can be used for the development of the active subspaces of functions with non-independent variables in next future.

\begin{appendices}  
\section{Proof of  Proposition \ref{prop:uniq}} \label{app:prop:uniq}
For continuous random variables and prescribed $\left(w_1, \ldots, w_{d-1} \right)$, the Rosenblatt transform of $\bo{\X}_{\sim j} |\X_j$ is unique and strictly increasing (\cite{rosenblatt52}). Therefore, the inverse of the Rosenblatt transform of $\bo{\X}_{\sim j} |\X_j$ is also unique (\cite{obrien75}), and  we can write  
$$
\bo{\X}_{\sim j} \stackrel{d}{=} r_{j}'\left(\X_{j}, \bo{U} \right)\, , 
$$
where $\bo{U} \sim \mathcal{U}\left(0,\, 1 \right)^{d-1}$. For the prescribed distribution of the $d-1$ innovation variables $\bo{Z}=\left(Z_{w_i} \sim F_{Z_{w_i}}, \, i=1, \ldots, d-1 \right)$, the above model becomes
   
$$
\bo{\X}_{\sim j} \stackrel{d}{=} r_{j}'\left(\X_{j}, U_1, \ldots, U_{d-1} \right) \stackrel{d}{=} r_{j}'\left(\X_{j}, F_{Z_{w_1}} \left(Z_{w_1}\right), \ldots, F_{Z_{w_{d-1}}} \left(Z_{w_{d-1}}\right) \right) = r_{j}\left(\X_{j}, \bo{Z} \right)\, , 
$$     
because      
$
Z_{w_i} \stackrel{d}{=} F_{Z{w_i}}^{-1}\left(U_{i} \right) 
\Longleftrightarrow U_{i} = F_{Z_{w_i}}\left(Z_{w_i} \right) 
$  
for continuous variable. Thus, Point (i) holds.\\   
For Point (ii), since $\bo{\X}_{\sim j} \stackrel{d}{=} r_{j}'\left(\X_{j}, \bo{U} \right)$ is the inverse of the Rosenblatt transform of $\bo{\X}_{\sim j} |\X_j$, we then have the unique inverse
$$
 \bo{U} = r_{j}^{'^{-1}}\left(\bo{\X}_{\sim j} |\X_j \right) \, , 
$$   
which yields to the unique inverse of the dependency model. Indeed, 
$$
(F_{\Z_1}(Z_1), \ldots, F_{\Z_{d-1}}(Z_1)) = r_{j}^{'^{-1}}\left(\bo{\X}_{\sim j} |\X_j \right)\Longrightarrow \, \bo{Z}  = r_{j}^{-1}\left(\bo{\X}_{\sim j} |\X_j \right) \, . 
$$       
   
\section{Proof of Lemma \ref{lem:unigrad}} \label{app:lem:unigrad}
Using the partial derivatives $\frac{\partial r_{\pi_{i,k}}}{\partial x_{\pi_{j,k}}}\left( \X_{\pi_{j,k}},\, \bo{Z}  \right)$ with $i=1, \ldots, d_k$ given by Equation (\ref{eq:fgrdc}) and the unique inverse of $\bo{\X}_{\sim \pi_{j,k}} \stackrel{}{=} r_{\pi_{j,k}}\left(\X_{\pi_{j,k}}, \bo{Z}_{k} \right)$ for any distribution of $\bo{Z}_{k}$ given by $\bo{Z}_{k} =r_{\pi_{j,k}}^{-1}\left(\bo{\X}_{\sim \pi_{j,k}} |\X_{\pi_{j,k}} \right)$ (see Proposition \ref{prop:uniq}), the result is immediate.   
      
\section{Proof of Theorem \ref{theo:agrad}} \label{app:theo:agrad}
Firstly, using the partial derivatives of $\bo{\X}_{\boldsymbol{\pi}_k}$ w.r.t $\X_{\pi_{j,k}}$ in (\ref{eq:unigrad}), that is,
$
J^{(\pi_{j,k})} =\left[\frac{\partial \X_{\pi_{1,k}}}{\partial x_{\pi_{j,k}}} \ldots \frac{\partial \X_{\pi_{d_k,k}}}{\partial x_{\pi_{j,k}}} \right]^T 
$
with 
$  
 \frac{\partial X_{\pi_{i,k}}}{\partial \x_{\pi_{j,k}}} = \frac{\partial r_{\pi_{i,k}}}{\partial x_{\pi_{j,k}}} =J^{(\pi_{j,k})}_{i}    
$
for any $\pi_{i,k} \in \boldsymbol{\pi}_k$, the reciprocal rule  allows for writing  
$$
\frac{\partial X_{\pi_{j,k}}}{\partial  \X_{\pi_{i,k}}} = \frac{1}{\frac{\partial \X_{\pi_{i,k}}}{\partial \X_{\pi_{j,k}}}} = \frac{1}{J^{(\pi_{j,k})}_{i} }\, .  
$$
Applying the chain rule yields   
$$
\frac{\partial X_{\pi_{i_1,k}}}{\partial  \X_{\pi_{i_2,k}}} = \frac{\partial \X_{\pi_{i_1,k}}}{\partial \X_{\pi_{j,k}}}\frac{\partial X_{\pi_{j,k}}}{\partial  \X_{\pi_{i_2,k}}}
 = \frac{J^{(\pi_{j,k})}_{i_1}}{J^{(\pi_{j,k})}_{i_2}} \, .       
$$     
Thus, the partial derivatives of $\bo{\X}_{\boldsymbol{\pi}_k}$ w.r.t. $\X_{\pi_{i,k}}$  are given by  
$$     
\frac{\partial \bo{\X}_{\boldsymbol{\pi}_k}}{\partial \X_{\pi_{i,k}}} := \left[\frac{J^{(\pi_{j,k})}_{1}}{J^{(\pi_{j,k})}_{i}}\; \ldots \; 
\frac{J^{(\pi_{j,k})}_{d_k}}{J^{(\pi_{j,k})}_{i}} \right]^T = \frac{J^{(\pi_{j,k})}}{J^{(\pi_{j,k})}_{i}}  \, ,     
$$          
and the actual Jacobian matrix of $\bo{\X}_{\boldsymbol{\pi}_k}$ (i.e., $\frac{\partial_a \bo{\X}_{\boldsymbol{\pi}_k}}{\partial \bo{\X}_{\boldsymbol{\pi}_k}}$) is given by 
$$    
J_{\boldsymbol{\pi}_k}^a := \left[ 
\begin{array}{cccc} 
\frac{J^{(\pi_{j,k})}}{J^{(\pi_{j,k})}_{1}}  & \ldots & \frac{J^{(\pi_{j,k})}}{J^{(\pi_{j,k})}_{d_k}}  \\
\end{array}
\right]        \, .      
$$ 
Secondly, under the organization of input variables (O) and for every pair $k_1,\, k_2 \in \{1, \ldots, K\}$ with $k_1 \neq k_2$, we have the following cross-Jacobian matrices: 
$$
\frac{\partial \bo{\X}_{\boldsymbol{\pi}_{k_1}}}{\partial \bo{x}_{\boldsymbol{\pi}_{k_2}}} =\mathsf{O}_{d_{k_1} \times d_{k_2}};
\qquad
J_{\boldsymbol{\pi}_1} := \frac{\partial \bo{\X}_{\boldsymbol{\pi}_1}}{\partial \bo{x}_{\boldsymbol{\pi}_1}} =\mathcal{I}_{d_1 \times d_1} \, .  
$$    
Therefore, the actual Jacobian matrix of $\bo{\X}$ is given by 
\begin{eqnarray}
J^a &:=& \left[  
\begin{array}{cccc}
\frac{\partial \bo{\X}_{\boldsymbol{\pi}_1}}{\partial \bo{x}_{\boldsymbol{\pi}_1}} & \frac{\partial \bo{\X}_{\boldsymbol{\pi}_1}}{\partial \bo{x}_{\boldsymbol{\pi}_2}}& \ldots & \frac{\partial \bo{\X}_{\boldsymbol{\pi}_1}}{\partial \bo{x}_{\boldsymbol{\pi}_K}}\\
\frac{\partial \bo{\X}_2}{\partial \bo{x}_{\boldsymbol{\pi}_1}} & \frac{\partial \bo{\X}_2}{\partial \bo{x}_{\boldsymbol{\pi}_2}}  &  \ldots & \frac{\partial \bo{\X}_2}{\partial \bo{x}_{\boldsymbol{\pi}_K}} \\
\vdots & \vdots  & \ddots & \vdots \\  
\frac{\partial \bo{\X}_{\boldsymbol{\pi}_K}}{\partial \bo{x}_{\boldsymbol{\pi}_1}} & \frac{\partial \bo{\X}_{\boldsymbol{\pi}_K}}{\partial \bo{x}_{\boldsymbol{\pi}_2}}& \ldots & \frac{\partial \bo{\X}_{\boldsymbol{\pi}_K}}{\partial \bo{x}_{\boldsymbol{\pi}_K}}\\
\end{array} 
 \right] = \left[         
\begin{array}{cccc}
\mathcal{I}_{d_1 \times d_1} & \mathsf{O}& \ldots & \mathsf{O} \\ 
\mathsf{O} & J_{\boldsymbol{\pi}_2}^a & \ldots &  \mathsf{O} \\
\vdots & \vdots  & \ddots & \vdots \\  
\mathsf{O} & \mathsf{O} & \ldots & J_{\boldsymbol{\pi}_K}^a\\
\end{array}           
 \right]  \nonumber      \, .  
\end{eqnarray} 
Finally, using the formal gradient of $\M$ in (\ref{eq:fgrdf}), that is,
$
\nabla \M := \left[\nabla\M_{\bo{x}_{\boldsymbol{\pi}_1}}^\T \; \, \nabla\M_{\bo{x}_{\boldsymbol{\pi}_2}}^\T \,  \ldots\ ,  \nabla\M_{\bo{x}_{\boldsymbol{\pi}_K}}^\T \right]^T
$
and bearing in mind the cyclic rule, we can write 
$$
\frac{\partial \M}{\partial \bo{x}_{\boldsymbol{\pi}_1}} = \left(\nabla \M^\T \frac{\partial \bo{\X}}{\partial \bo{x}_{\boldsymbol{\pi}_1}}\right)^\T =\nabla \M_{\bo{x}_{\boldsymbol{\pi}_1}} ,
\qquad
\frac{\partial \M}{\partial \bo{x}_{\boldsymbol{\pi}_k}} = \left(\nabla \M^\T \frac{\partial \bo{\X}}{\partial \bo{x}_{\boldsymbol{\pi}_k}}\right)^\T = \left(\nabla \M_{\bo{x}_{\boldsymbol{\pi}_k}}^\T J_{\boldsymbol{\pi}_k}^a \right)^\T   \, ,  
$$          
and the actual partial derivatives of $\M$ are then given by   
$   
\frac{\partial_a \M}{\partial \bo{x}}  = J^a(\bo{x})^\T \nabla \M(\bo{x})
$.    

\section{Proof of Theorem \ref{theo:dgrad}} \label{app:theo:dgrad}
For Point (i), building the $d_k$ dependency functions for every explanatory input $\X_{\pi_{i,k}}$ with $i=1, \ldots, d_k$, the partial derivatives of $\bo{\X}_{\boldsymbol{\pi}_k}$ w.r.t $\X_{\pi_{i,k}}$ evaluated at $\bo{x}_{\boldsymbol{\pi}_k}$ is given by (see Equation (\ref{eq:unigrad}))                                                                            
$$
J^{(\pi_{i,k})}(\bo{x}_{\boldsymbol{\pi}_k}) = J^{(\pi_{i,k})}\left((x_{\pi_{i,k}},\, r_{\pi_{i,k}}^{-1}\left(\bo{x}_{\pi_{i,k}}|x_{\pi_{i,k}} \right) \right); \quad  i=1, \ldots, d_k \, .
$$ 
Points (ii)-(iii) are similar to those of Theorem \ref{theo:agrad} using the dependent Jacobian matrix given by Point (i). 

\section{Proof of Theorem \ref{theo:dhess}} \label{app:theo:dhess}
First, using Equation (\ref{eq:dgradf}) given by
$
\frac{\partial \M}{\partial \bo{x}}(\bo{x})  :=  J^d \left(\bo{x} \right)^\T  \nabla \M(\bo{x})    
$, 
we can extract
$$ 
\frac{\partial \M}{\partial \bo{x}_{\boldsymbol{\pi}_1}} (\bo{x}) =\nabla \M_{\bo{x}_{\boldsymbol{\pi}_1}}(\bo{x}); 
\qquad   \quad 
\frac{\partial \M}{\partial \bo{x}_{\boldsymbol{\pi}_k}} (\bo{x}) = J_{\boldsymbol{\pi}_k}^d\left(\bo{x}_{\boldsymbol{\pi}_k} \right)^\T \nabla \M_{\bo{x}_{\boldsymbol{\pi}_k}}(\bo{x}) \, , 
$$ 
By applying the vector-by-vector derivatives of $\frac{\partial \M}{\partial \bo{x}_{\boldsymbol{\pi}_1}} (\bo{x})$ w.r.t. $\bo{x}_{\boldsymbol{\pi}_1}$ and $\bo{x}_{\boldsymbol{\pi}_k}$, we have  
$$
\frac{\partial^2 \M}{\partial^2 \bo{x}_{\boldsymbol{\pi}_1}}(\bo{x}) = 
\frac{\partial \left[\nabla \M_{\bo{x}_{\boldsymbol{\pi}_1}} \right]}{\partial \bo{x}_{\boldsymbol{\pi}_1}}(\bo{x}) = H_{\boldsymbol{\pi}_1}(\bo{x}); 
\qquad     
\frac{\partial^2 \M(\bo{x})}{\partial \bo{x}_{\boldsymbol{\pi}_1} \partial \bo{x}_{\boldsymbol{\pi}_k}} = 
\frac{\partial \left[\nabla \M_{\bo{x}_{\boldsymbol{\pi}_1}}(\bo{x}) \right]}{\partial \bo{x}} \frac{\partial \bo{x}}{\partial \bo{x}_{\boldsymbol{\pi}_k}}  
= H_{\boldsymbol{\pi}_1,\boldsymbol{\pi}_k}(\bo{x}) J_{\boldsymbol{\pi}_k}^d \,  ,
$$         
as $\bo{\X}_{\boldsymbol{\pi}_1}$ is a vector of independent variables and 
$ 
\frac{\partial \bo{\X}}{\partial \bo{x}_{\boldsymbol{\pi}_k}} = \left[ 
\begin{array}{c}
\mathsf{O}_{d_1 \times d_k} \\ 
\vdots \\
J_{\boldsymbol{\pi}_k}^d \\ 
\vdots \\
 \mathsf{O}_{d_K \times d_k} \\ 
\end{array}
\right] 
$, bearing in mind the dependent Jacobian matrix provided in Equation (\ref{eq:djacall}). \\
In the same sense, the derivatives of $\frac{\partial \M}{\partial \bo{x}_{\boldsymbol{\pi}_k}}(\bo{x})$ w.r.t. $\bo{x}_{\boldsymbol{\pi}_1}$ and $\bo{x}_{\boldsymbol{\pi}_\ell}$ with $\ell \neq k$ are 
   
\begin{eqnarray}  
\frac{\partial^2 \M}{\partial \bo{x}_{\boldsymbol{\pi}_k} \partial \bo{x}_{\boldsymbol{\pi}_1}} &=& 
\frac{\partial \left[J_{\boldsymbol{\pi}_k}^d\left(\bo{x}_{\boldsymbol{\pi}_k} \right)^\T \nabla \M_{\bo{x}_{\boldsymbol{\pi}_k}}(\bo{x}) \right]}{\partial \bo{x}_{\boldsymbol{\pi}_1}}  
= J_{\boldsymbol{\pi}_k}^d\left(\bo{x}_{\boldsymbol{\pi}_k} \right)^\T \frac{\partial \left[\nabla \M_{\bo{x}_{\boldsymbol{\pi}_k}}(\bo{x}) \right]}{\partial \bo{x}_{\boldsymbol{\pi}_1}}  \nonumber \\
&=& J_{\boldsymbol{\pi}_k}^d\left(\bo{x}_{\boldsymbol{\pi}_k} \right)^\T \frac{\partial \left[\nabla \M_{\bo{x}_{\boldsymbol{\pi}_k}}(\bo{x}) \right]}{\partial \bo{x}}  \frac{\partial  \bo{x}}{\partial \bo{x}_{\boldsymbol{\pi}_1}} 
= J_{\boldsymbol{\pi}_k}^d\left(\bo{x}_{\boldsymbol{\pi}_k} \right)^\T H_{\boldsymbol{\pi}_k, \boldsymbol{\pi}_1}(\bo{x}) \,  ,  \nonumber          
\end{eqnarray}      
\begin{eqnarray} 
\frac{\partial^2 \M}{\partial \bo{x}_{\boldsymbol{\pi}_k} \partial \bo{x}_{\boldsymbol{\pi}_\ell}} &=& 
\frac{\partial \left[J_{\boldsymbol{\pi}_k}^d\left(\bo{x}_{\boldsymbol{\pi}_k} \right)^\T \nabla \M_{\bo{x}_{\boldsymbol{\pi}_k}}(\bo{x}) \right]}{\partial \bo{x}_{\boldsymbol{\pi}_\ell}} =
  J_{\boldsymbol{\pi}_k}^d\left(\bo{x}_{\boldsymbol{\pi}_k} \right)^\T \frac{\partial \left[\nabla \M_{\bo{x}_{\boldsymbol{\pi}_k}}(\bo{x}) \right]}{\partial \bo{x}_{\boldsymbol{\pi}_\ell}} \nonumber \\ 
&=& J_{\boldsymbol{\pi}_k}^d\left(\bo{x}_{\boldsymbol{\pi}_k} \right)^\T \frac{\partial \left[\nabla \M_{\bo{x}_{\boldsymbol{\pi}_k}}(\bo{x}) \right]}{\partial \bo{x}} \frac{\partial \bo{x}}{\partial \bo{x}_{\boldsymbol{\pi}_\ell}}  
	= J_{\boldsymbol{\pi}_k}^d\left(\bo{x}_{\boldsymbol{\pi}_k} \right)^\T  
H_{\boldsymbol{\pi}_k, \boldsymbol{\pi}_\ell}(\bo{x}) J_{\bo{x}_{\boldsymbol{\pi}_\ell}}^d (\bo{x}_{\boldsymbol{\pi}_\ell})\,  . 
\nonumber      
\end{eqnarray}   
Finally, we have to derive the quantity  
$
\frac{\partial^2 \M}{\partial^2 \bo{x}_{\boldsymbol{\pi}_k}} = 
\frac{\partial \left[J_{\boldsymbol{\pi}_k}^d\left(\bo{x}_{\boldsymbol{\pi}_k} \right)^\T \nabla \M_{\bo{x}_{\boldsymbol{\pi}_k}}(\bo{x}) \right]}{\partial \bo{x}_{\boldsymbol{\pi}_k}} 
$.  
For each $\pi_{\ell, k} \in \boldsymbol{\pi}_k$,  we can write
      
\begin{eqnarray}   
& & \frac{\partial^2 \M}{\partial \bo{x}_{\boldsymbol{\pi}_k} \partial x_{\pi_{\ell, k}}} = 
\frac{\partial \left[J_{\boldsymbol{\pi}_k}^d\left(\bo{x}_{\boldsymbol{\pi}_k} \right)^\T \nabla \M_{\bo{x}_{\boldsymbol{\pi}_k}}(\bo{x}) \right]}{\partial x_{\pi_{\ell, k}}} \nonumber \\  
&=& \frac{\partial \left[J_{\boldsymbol{\pi}_k}^d\left(\bo{x}_{\boldsymbol{\pi}_k} \right)^\T \right]}{\partial x_{\pi_{\ell, k}}}  \nabla \M_{\bo{x}_{\boldsymbol{\pi}_k}}(\bo{x})  +   J_{\boldsymbol{\pi}_k}^d\left(\bo{x}_{\boldsymbol{\pi}_k} \right)^\T 
\frac{\partial \left[\nabla \M_{\bo{x}_{\boldsymbol{\pi}_k}}(\bo{x}) \right]}{\partial x_{\pi_{\ell, k}}} \nonumber \\
&=&  \left[\frac{\partial J^{(\pi_{1,k})}}{\partial x_{\pi_{\ell, k}}} \ldots   \frac{\partial J^{(\pi_{d_k,k})}}{\partial x_{\pi_{\ell, k}}} \right]^\T \nabla \M_{\bo{x}_{\boldsymbol{\pi}_k}}(\bo{x}) + J_{\boldsymbol{\pi}_k}^d\left(\bo{x}_{\boldsymbol{\pi}_k} \right)^\T    
\frac{\partial \left[\nabla \M_{\bo{x}_{\boldsymbol{\pi}_k}}(\bo{x}) \right]}{\partial \bo{x}} 
\frac{\partial \bo{x}}{{\partial x_{\pi_{\ell, k}}}}
\nonumber \\       
&=& \left[J^{(\pi_{\ell, k})}_1  \frac{\partial J^{(\pi_{1,k})}}{\partial x_{\pi_{1, k}}} \ldots J^{(\pi_{\ell, k})}_{d_k}  \frac{\partial J^{(\pi_{d_k,k})}}{\partial x_{\pi_{d_k, k}}} \right]^\T \nabla \M_{\bo{x}_{\boldsymbol{\pi}_k}}(\bo{x}) + J_{\boldsymbol{\pi}_k}^d\left(\bo{x}_{\boldsymbol{\pi}_k} \right)^\T   H_{\boldsymbol{\pi}_k}(\bo{x}) J^{(\pi_{\ell, k})} (\bo{x}_{\boldsymbol{\pi}_k})      \nonumber \, ,    
\end{eqnarray}   
because for all $i \in\{1,\ldots, d_k\}$, we can write (thanks to the chain rule)  
$$ 
\frac{\partial J^{(\pi_{i,k})}}{\partial x_{\pi_{\ell, k}}} :=\left[\frac{\partial^2 \X_{\pi_{1, k}}}{\partial^2 x_{\pi_{i,k}}} \frac{\partial x_{\pi_{i,k}}}{\partial x_{\pi_{\ell, k}}}     
\; \ldots \; \frac{\partial^2 \X_{\pi_{d_k, k}}}{\partial^2 x_{\pi_{i,k}}} 
\frac{\partial x_{\pi_{i,k}}}{\partial x_{\pi_{\ell, k}}} \right]^T =   
\frac{\partial x_{\pi_{i,k}}}{\partial x_{\pi_{\ell, k}}}  \frac{\partial J^{(\pi_{i,k})}}{\partial x_{\pi_{i,k}}} 
=  J^{(\pi_{\ell, k})}_i  \frac{\partial J^{(\pi_{i,k})}}{\partial x_{\pi_{i,k}}} \, .  
$$  
Re-organizing the first element of the right-hand terms of the above equation yields
\begin{eqnarray} 
\frac{\partial^2 \M}{\partial \bo{x}_{\boldsymbol{\pi}_k} \partial x_{\pi_{\ell, k}}}       
&=&  diag\left( \left[\frac{\partial J^{(\pi_{1,k})}}{\partial x_{\pi_{1, k}}} \ldots  \frac{\partial J^{(\pi_{d_k,k})}}{\partial x_{\pi_{d_k, k}}} \right]^\T\nabla \M_{\bo{x}_{\boldsymbol{\pi}_k}}(\bo{x}) \right) J^{(\pi_{\ell, k})} + J_{\boldsymbol{\pi}_k}^d\left(\bo{x}_{\boldsymbol{\pi}_k} \right)^\T   H_{\boldsymbol{\pi}_k}(\bo{x}) J^{(\pi_{\ell, k})}    
\nonumber \, .  
\end{eqnarray}  
By running $\ell =1, \ldots, d_k$, we obtain the result. 
        
\end{appendices}

                        

                                           
\end{document}